\documentclass[11pt]{article}
\usepackage{amsmath}

\usepackage{relsize}
\usepackage{color}
\usepackage{appendix}
\title{{\bf On a double integral of a product of Legendre polynomials}}
\author{ G. Vaman\\Institute of Atomic Physics, P. O. Box MG-6, Bucharest, Romania\\ (vaman@ifin.nipne.ro)}
\begin{document}
\maketitle
\begin{abstract}
\noindent
We calculate a double integral over a product of Legendre polynomials multiplied by a binomial raised to a power.
\end{abstract}

During the calculation of the electromagnetic self-force of a uniformly charged spherical ball, we encountered the integral
\begin{align}\label{1}
I=\int_0^{\pi}d \theta \sin \theta \int_0^{\pi}d \theta' \sin \theta' (\cos \theta - \cos \theta')^{2n} P_l(\cos \theta) P_l(\cos \theta'),
\end{align}
where $n$ and $l$ are integer positive numbers and $P_l$ is a Legendre polynomial of order $l$. As far as we know, this integral was not calculated in closed form anywhere in the literature. We calculate it here.

After changing the variables $\cos \theta \rightarrow x$, $\cos \theta' \rightarrow y$ this integral becomes
\begin{align}\label{2}
I=\int_{-1}^1 dx \int_{-1}^1 dy \;(x-y)^{2n} P_l(x) P_l(y).
\end{align}
We first perform the integral $\int_{-1}^1 dx (x-y)^{2n} P_l(x).$
For this, we can use {\bf 7.228} and {\bf 8.703} from \cite{gra}. Combining these two equations that read
\begin{align}
\frac{1}{2}\Gamma(1+\mu) \int_{-1}^1 P_l(x) (z-x)^{-\mu-1}= (z^2-1)^{- \frac{\mu}{2}} e^{-i \pi \mu} Q_l^{\mu}(z), l=0,1,\dots, |arg(z-1)|<\pi,
\end{align}
\begin{align}
Q_{\nu}^{\mu}(x)= \frac{e^{i\pi \mu} \Gamma(\nu+\mu+1) \Gamma\left(\frac{1}{2}\right)}{2^{\nu+1} \Gamma\left( \nu + \frac{3}{2}\right)} (x^2-1)^{\frac{\mu}{2}} x^{-\nu-\mu-1} {}_2F_{1}\left(\frac{\nu+ \mu +2}{2}, \frac{\nu+\mu+1}{2}; \nu+\frac{3}{2}; \frac{1}{x^2}\right),
\end{align}
after using \[ \frac{\Gamma(l+\mu+1)}{\Gamma(\mu+1)} =(\mu +1)_l\; \text{and} \; \Gamma\left( l+\frac{3}{2}\right)= \frac{\sqrt{\pi} (l+1)_{l+1}}{2^{2l+1}},\]
one obtains for $\mu=-1-2n$
\begin{align}
\int_{-1}^1 dx \, (x-y)^{2n}P_l(x)  = \frac{(-2n)_l \,2^{l+1} y^{2n-l}}{(l+1)_{l+1}} {}_2F_1 \left( \frac{l}{2}-n+\frac{1}{2}, \frac{l}{2}-n; l+\frac{3}{2}; \frac{1}{y^2} \right).
\end{align}
The same result given in Eq. (5) can be obtained by putting $a=1$, $m=0$ and $p=-2n$ in Eq. {\bf 2.17.4(5)} from \cite{pru} 
\begin{align}
\int_{-a}^a dx \,\frac{(a^2-x^2)^{\frac{m}{2}}}{(x-y)^p} P_l^m \left( \frac{x}{a} \right)= \frac{2 (-1)^{m-1}(l+m)!}{(p-1)! (l-m)!} (y^2-a^2)^{\frac{m-p+1}{2}} Q_l^{p-m-1} \left( \frac{y}{a}\right),
\end{align}
although this equation is given in \cite{pru} as being valid only for $p=0,1, \dots$.

The same result given in Eq.(5) can be obtained by direct calculation, by using the Rodrigues formula for Legendre polynomials \cite{rai}
\begin{align}
P_l(x)= \sum_{k=0}^{[l/2]} \frac{(-1)^k \left( \frac{1}{2}\right)_{l-k} (2x)^{l-2k}}{k!(l-2k)!}
\end{align}
and the binomial expansion for $(x-y)^{2n}$
\begin{align}
(x-y)^{2n}= \sum_{p=0}^{2n} \frac{(-1)^p (2n)!}{p!(2n-p)!} x^p y^{2n-p},
\end{align}
 and integrating the resulting double sum term by term. For $l$ odd, after noting that the term by term integration gives non-zero result only for $p$ odd and changing the summation index $p \rightarrow 2p$, one obtains
\begin{align}
\int_{-1}^1 dx \,(x-y)^{2n} P_l(x)= -2^{l+1} \sum_{k=0}^{[l/2]} \sum_{p=0}^{b-1} \frac{(-1)^k \left( \frac{1}{2} \right)_{l-k}(2n)! y^{2n-2p-1}}{k! (l-2k)!2^{2k}(2p+1)!(2n-2p-1)!}\nonumber \\
\cdot \frac{1}{(l-2k+2p+2)}.
\end{align}
Writing all the factorials in terms of Pochhammer symbols, the above summation over $k$ can be done as follows 
\begin{align}
&\sum_{k=0}^{\frac{l-1}{2}} \frac{(-1)^k \left( \frac{1}{2}\right)_{l-k}}{k!(l-2k)!2^{2k}(l+2p+2-2k)} \nonumber \\
&= \frac{\left( \frac{1}{2}\right)_l}{\Gamma(1+l) (l+p+2)} {}_3F_2 \left( -\frac{l}{2}+\frac{1}{2}, -\frac{l}{2}, -\frac{l}{2}-p-1; \frac{1}{2}-l, -\frac{l}{2}-p,1\right)\nonumber \\
&\frac{\left(\frac{1}{2}\right)_l \left( \frac{1-l}{2} \right)_{\frac{l-1}{2}} (-p)_{\frac{l-1}{2}}}{l! \,(l+2p+2) \left(\frac{1}{2}-l\right)_{\frac{l-1}{2}} \left(-\frac{l}{2}-p\right)_{\frac{l-1}{2}}},
\end{align}
where, when we passed from the second to the third line of the above equation, we used equation {\bf 7.4.4 (81)} from \cite{pru}. We note that, because of the Pochhammer symbol $(-p)_{\frac{l-1}{2}}$, the r.h.s. of Eq. (10) is different from zero only for $p \ge \frac{l-1}{2}$. Introducing Eq. (10) in Eq. (9) and changing the summation index $p\rightarrow i, \, p-\frac{l-1}{2}=i$, the resulting summation over i can be done immediately and one obtains again the result of Eq. (5). The case $l$ even can be considered similarly.

Returning now to Eq. (2), after using Eq. (5) one obtains
\begin{align}\label{7}
I=\frac{(-2n)_l \;2^{l+1}}{(l+1)_{l+1}} \int_{-1}^1 dy\; y^{2n-l} P_l(y)\; {}_2F_1 \left(\frac{l}{2}-n, \frac{l+1}{2}-n, l+\frac{3}{2};\frac{1}{y^2} \right).
\end{align}
Note that, for $l \le 2n$, the hypergeometric function in Eq.(\ref{7}) is, in fact, a finite series, because $l/2-n$ and $(l+1)/2-n$ are negative integers when $l$ is even and odd respectively. Again, we consider the cases $l$ even and $l$ odd separately.

For $l$ even, using the definition of the Gauss hypergeometic function, we have
\begin{align}\label{8}
{}_2F_1 \left(\frac{l}{2}-n, \frac{l+1}{2}-n, l+\frac{3}{2};\frac{1}{y^2} \right)= \sum_{k=0}^{n-\frac{l}{2}} \frac{ \left( \frac{l}{2}-n \right)_k \left(\frac{l+1}{2}-n \right)_k}{k! \left( l+\frac{3}{2} \right)_k} \frac{1}{y^{2k}}.
\end{align}
From Eqs. (\ref{7}), (\ref{8}), one obtains
\begin{align}\label{9}
I= \frac{2 (-2n)_l2^{l+1}}{(l+1)_{l+1}} \sum_{k=0}^{n-\frac{l}{2}} \frac{ \left( \frac{l}{2}-n \right)_k \left(\frac{l+1}{2}-n \right)_k}{k! \left( l+\frac{3}{2} \right)_k} \int_0^1 dy \; y^{2n-2k-l} P_l(y),
\end{align}
where we used the fact that the integrand in the r.h.s. of Eq. (\ref{9}) is an even function, because $P_l(-y)=(-1)^l P_l(y)$ \cite{gra}. The integral in Eq. (\ref{9}) can be performed using ({\bf 2.17.1(4)} from \cite{pru} or {\bf 7.126(2)} from \cite{gra}  )
\begin{align}\label{p2}
\int_0^1dx \; x^{\sigma} P_{\nu}(x)= \frac{\sqrt{\pi} \;2^{-\sigma-1} \Gamma(1+\sigma)}{\Gamma \left(1+\frac{\sigma-\nu}{2} \right) \Gamma \left( \frac{\sigma+\nu+3}{2} \right)}.
\end{align}
One obtains
\begin{align}\label{11}
I=\frac{\sqrt{\pi} \;(-2n)_l \;2^{2l-2n+1}}{(l+1)_{l+1}} \sum_{k=0}^{n-\frac{l}{2}} \frac{2^{2k} \left( \frac{l}{2}-n \right)_k \left( \frac{l+1}{2}-n \right)_k \Gamma(2n-2k-l+1)}{k! \left( l+\frac{3}{2} \right)_k \Gamma(1+n-k-l) \Gamma \left( n-k+\frac{3}{2} \right)}.
\end{align}
Using the definition of the Pochhammer symbol \cite{pru} \[ (a)_k= \frac{\Gamma(a+k)}{\Gamma(a)}= (-1)^k \frac{\Gamma(1-a)}{\Gamma(1-a-k)}, \]
and \cite{pru}
\[(a)_{2k}= \left( \frac{a}{2} \right)_k \left( \frac{a+1}{2} \right)_k 2^{2k}, \]
we write the Gamma functions in Eq. (\ref{11}) as follows
\begin{align}
& \Gamma(2n-2k-l+1)= \frac{\Gamma(2n-l+1)}{\left( \frac{l}{2}-n \right)_k \left( \frac{l}{2}-n+\frac{1}{2}\right)_k2^{2k}}, \nonumber\\
&\Gamma(1+n-l-k)= (-1)^k\frac{\Gamma(1+n-l)}{(l-n)_k},\nonumber \\
&\Gamma \left(n+\frac{3}{2}-k \right)= (-1)^k \frac{\Gamma \left( n+\frac{3}{2}\right)}{\left( -n-\frac{1}{2}\right)_k}.
\end{align}
Introducing (16) in (15), one obtains
\begin{align}
I=\frac{\sqrt{\pi} \;(-2n)_l\; 2^{2l-2n+1}(2n-l)!}{ (l+1)_{l+1} (n-l)! \left( n+\frac{1}{2} \right)!} {}_2F_1 \left(-n-\frac{1}{2}, l-n; l+\frac{3}{2}; 1 \right).
\end{align}
But the Gauss hypergeometric function of unit argument can be written as ({\bf 7.3.5(2)} \cite{pru})
\begin{align}
{}_2F_1(a,b;c;1)=\frac{\Gamma(c) \Gamma(c-a-b)}{\Gamma(c-a) \Gamma(c-b)},
\end{align}
and we obtain
\begin{align}\label{17}
I=\frac{ \sqrt{\pi}\; (-2n)_l \;2^{2l-2n+1} (2n-l)! \left(l+\frac{1}{2}\right)! (2n+1)!}{ (l+1)_{l+1} (n-l)! \left(n+\frac{1}{2}\right)! (l+n+1)! \left(n+ \frac{1}{2} \right)!},
\end{align}
where we use the notation $\Gamma(z)=(z-1)!$ both for integer and noninteger $z$. Using the definition of the Pochhammer symbol and \cite{pru}
\begin{align}
\frac{\Gamma(2z)}{\Gamma(z)}=\frac{2^{2z-1}}{\sqrt{\pi}} \Gamma \left( z+\frac{1}{2} \right),
\end{align}
we can write 
\begin{align}\label{19}
&(-2n)_l= (-1)^l \frac{\Gamma(1+2n)}{\Gamma(1+2n-l)},\nonumber \\
& (l+1)_{l+1}=\frac{2^{2l+1}}{\sqrt{\pi}} \Gamma \left(l+\frac{3}{2} \right),\nonumber \\
& \frac{(2n)!}{\left( n+\frac{1}{2}\right)!}=\frac{2^{2n+1}\Gamma(n+1)}{\sqrt{\pi} (2n+1)},\\
& \frac{(2n+1)!}{\left( n+\frac{1}{2}\right)!}= \frac{2^{2n+1}\Gamma(n+1)}{\sqrt{\pi} }. \nonumber
\end{align}
Introducing (\ref{19}) in (\ref{17}), one obtains for $l$ even 
\begin{align}
I= \frac{(-1)^l 2^{2n+2} (n!)^2}{(2n+1) (n-l)! (n+l+1)!}.
\end{align}
Note that, because of $(n-l)!$ from the denominator, this result is different from zero only for $n \ge l$. A similar calculation can be done for $l$ odd, and one obtains the same result. So, our final result for the integral (\ref{1}) is
\begin{equation}
I= \left\lbrace\begin{array}{c}
               \frac{(-1)^l 2^{2n+2} (n!)^2}{(2n+1) (n-l)! (n+l+1)!}, n \ge l\\
               0, n<l
               \end{array}
                       \right..
\end{equation}

\end{document}